

\nopagenumbers
\parindent= 15pt

\hsize=14cm
\vsize=21cm
\hoffset=0.9cm
\voffset=1cm

\input amssym.def
\input amssym.tex

 at6.5pt
\font\srm=cmr8

\font\csc=cmcsc10
\font\title=cmr12 at 14pt

\font\teneusm=eusm10    
\font\seveneusm=eusm7  
\font\fiveeusm=eusm5    
\newfam\eusmfam

\textfont\eusmfam=\teneusm 
\scriptfont\eusmfam=\seveneusm
\scriptscriptfont\eusmfam=\fiveeusm

\def\Re{{\rm Re}\,}
\def\Im{{\rm Im}\,}

\def\scr#1{{\scriptstyle{#1}}}

\def\rightheadline{\hfil{\srm On Some Improvements 
of the Brun--Titchmarsh Theorem. IV}\hfil\tenrm\folio}
\def\leftheadline{\tenrm\folio\hfil{\srm Y. Motohashi}\hfil}
\def\emptyheadline{}
\headline{\ifnum\pageno=1 \emptyheadline\else
\ifodd\pageno \rightheadline \else \leftheadline\fi\fi}

\def\firstpage{\hss{\vbox to 1.5cm
{\vfil\hbox{\rm\folio}}}\hss}
\def\emptyfootline{\hfil}
\footline{\ifnum\pageno=1\firstpage\else
\emptyfootline\fi}

\centerline{\title On Some Improvements of the Brun--Titchmarsh
Theorem. IV}
\vskip 0.5cm
\centerline{\csc By Yoichi Motohashi}
\vskip 0.7cm
\noindent
The aim of the present note is to give a version of large sieve
extensions of the Brun--Titchmarsh theorem. This is in fact a
rework of our old file left unpublished since early 1980's 
which we originally intended to include into our 
Tata lecture notes [12].
We publish it here in a slightly reworked form, as it might have now
some fresh interest in view of Maynard's recent publication [4].
\bigskip
\noindent
{\bf 1.} Let 
$$
\pi(x;k,\ell)=\sum_{\scr{p\le x}\atop\scr{p\equiv \ell\bmod k}}
1,\quad (k,\ell)=1, \eqno(1.1)
$$
where $p$ denotes a generic prime, and let $\chi$ be a Dirichlet
character. Then [12, Theorem 13] asserts, among other things, that
we have, uniformly for $kQ^2\le x^{9/20-\varepsilon}$,
$$
\sum_{\scr{q\le Q}\atop\scr{(q,k)=1}}\mathop{{\sum}^*}_{\chi\bmod q}
\bigg|\sum_{\scr{p\le x}\atop\scr{p\equiv\ell\bmod k}}
\chi(p)\bigg|^2
\le{(2+o(1))x\over\varphi(k)
\log(x/(kQ^2)^{3/8})}\pi(x;k,\ell),
\eqno(1.2)
$$
provided $x$ is larger than a constant which is effectively
computable for each $\varepsilon>0$,
where $\varphi$ is the Euler totient function and
$\sum^*$ stands for the restriction of the sum to primitive characters. 
In particular we have
$$
\pi(x;k,\ell)\le{(2+o(1))x\over\varphi(k)\log(x/k^{3/8})},
\quad k\le x^{9/20-\varepsilon},\eqno(1.3)
$$
which surpasses partly the famed bound
$$
\pi(x;k,\ell)\le {2x\over\varphi(k)\log(x/k)},\quad k<x,\eqno(1.4)
$$
due to Montgomery and Vaughan [5].
In contrast to this, Maynard [4] asserts in essence that
$$
\pi(x;k,\ell)\le{2x\over\varphi(k)\log x},\quad k\le x^{1/8},
\eqno(1.5)
$$
provided $x$ is larger than an effectively 
computable constant. He gives also a lower bound, though
we skip it in order to make our presentation simple; 
for the same reason, we also skip mentioning 
former improvements upon $(1.4)$ other than $(1.3)$.
The bound $(1.5)$ has been known as a kind of folklore 
among specialists, but with less precision
about the range of moduli. 
\medskip
We shall refine $(1.2)$ by
\medskip
\noindent
{\bf Theorem.} {\it There exists an effectively computable
constant $\omega$ such that we have, uniformly for
$kQ^2\le x^\omega$,
$$
\sum_{\scr{q\le Q}\atop\scr{(q,k)=1}}
\mathop{{\sum}^*}_{\chi\bmod q}
\bigg|\sum_{\scr{p\le x}\atop\scr{p\equiv\ell\bmod k}}
\chi(p)\bigg|^2\le{2x\over\varphi(k)\log x}\pi(x;k,\ell).
\eqno(1.6)
$$
}
\medskip
\noindent
Obviously this contains $(1.5)$ but for $k\le x^\omega$. It 
remains thus to find a good lower bound 
for $\omega$. We are certain that
Maynard's argument will extend to 
the direction indicated by $(1.6)$ and yield $(1.5)$ as a particular 
instance, since the basic structure of 
his argument is essentially the same as ours that is
developed in [11][12] in detail, although the intricate
part of [4] corresponding to the numerical
precision should be overhauled accordingly. 
Further, we add that it is very possible
to prove a short interval version of $(1.6)$.
\medskip
\noindent
{\csc Remark:} The web edition of [12] contains some 
obvious misprints; for example the statement
there corresponding to $(1.2)$ lacks the necessary
restriction to primitive characters. We are going to provide
Tata IFR with corrections.
\bigskip
\noindent
{\bf 2.} Our theorem is in fact a simple 
consequence of our version [11][12, Theorem 17] of the
Linnik--Fogels--Gallagher prime number theorem; 
thus we need first
to introduce a notion concerning zeros of Dirichlet $L$-functions
$L(s,\chi)$. Details can be found in [16], for instance.
\medskip
We consider the set $Z_T$ of all non-trivial zeros 
in the region $|\Im s|\le T$ of the function
$\prod_{q\le T}\prod^*_{\chi\bmod q}L(s,\chi)$, i.e., with
$\chi$ being primitive. Then we have that there exists an
effectively computable absolute constant $a_0>0$
such that
$$
\hbox{$\displaystyle
\max_{\rho\in Z_T}\Re\rho\le 1-{a_0\over\log T}$,
excepting a possibly existing zero $\beta$.
}\eqno(2.1)
$$
If $\beta$ exists, it is real and simple, and we describe both itself and
the relevant unique primitive character as $T$-exceptional. 
We put
$$
\Delta_T=\cases{\hfil1&if $\beta$ does not exist,\cr 
(1-\beta)\log T & if $\beta$ exists.}\eqno(2.2)
$$
It is generally believed that $\beta$ 
does not exist. A way to confirm this 
is to improve appropriately the Brun--Titchmarsh theorem,
which remains, however, to be one of the most difficult problems in
analytic number theory; see [9][12, \S4.3].
\medskip
We put
$$
\psi(x,\chi)=\sum_{n\le x}\chi(n)\Lambda(n),\eqno(2.3)
$$
with the von Mangold function $\Lambda$;
and let
$$
\tilde{\psi}(x,\chi)=\cases{\psi(x,\chi)-x& if $\chi$ is principal,
\cr \psi(x,\chi)+x^\beta/\beta& if $\chi$ is $T$-exceptional,
\cr\psi(x,\chi)& otherwise.}\eqno(2.4)
$$
Then, [12, Theorem 17] asserts that there exist effectively
computable absolute constants $a_1, a_2, a_3>0$ such that
$$
\sum_{q\le T}\mathop{{\sum}^*}_{\chi\bmod q}
|\tilde{\psi}(x,\chi)|\le a_1 x\Delta_T
\exp\big(-a_2\log x/\log T\big),\eqno(2.5)
$$
provided $T^{a_3}\le x\le\exp((\log T)^2)$.
\medskip
In order to prove our theorem, we note that
$$
\eqalignno{
\sum_{\scr{p\le x}\atop\scr{p\equiv\ell\bmod k}}
\chi(n)\Lambda(n)&={1\over\varphi(k)}\sum_{\xi\bmod k}
\bar\xi(\ell)\psi(x,\xi\chi)\cr
&={1\over\varphi(k)}\sum_{\xi\bmod k}\bar{\xi}(\ell)
\psi(x,\xi^\sharp\chi)+O(\nu(k)\log x),&(2.6)
}
$$
where $\xi^\sharp$ is the primitive character inducing 
the Dirichlet character $\xi$, and
$\nu(k)$ the number of distinct prime factors of $k$.
Here $\xi^\sharp\chi$ stands for a unique primitive
character whose conductor is not larger than
$kQ$. We have thus
$$
\eqalignno{
&\sum_{\scr{q\le Q}\atop\scr{(q,k)=1}}
\mathop{{\sum}^*}_{\chi\bmod q}
\bigg|\sum_{\scr{p\le x}\atop\scr{p\equiv\ell\bmod k}}
\chi(n)\Lambda(n)\bigg|\cr
&\le{1\over\varphi(k)}\left\{x+x^\beta/\beta
+E(x,kQ)\right\}+O(\nu(k)Q^2\log x),&(2.7)
}
$$
where $\beta$ is the $kQ$-exceptional zero if exists;
and $E(x,kQ)$ is the left side of $(2.5)$ for
$T=kQ$. If $\beta$ exists, then $(2.5)$ implies that
$$
\eqalignno{
x^{-1}\left(x^\beta/\beta+E(x,kQ)\right)\le&\,
\exp\big(-\Delta_T\log x/\log T\big)/
(1-\Delta_T/\log T)\cr
+&\,a_1\Delta_T\exp\big(-a_2\log x/\log T\big),&(2.8)
}
$$
provided $T^{a_3}\le x\le((\log T)^2)$. The right side is
$$
\eqalignno{
&\le\exp\big(-a_3\Delta_T)/(1-\Delta_T/\log T)
+a_1\Delta_T\exp(-a_2a_3)\cr
&\le\exp(-a_3\Delta_T)
+\Delta_T\big(1/(2\log T)+a_1\exp(-a_2a_3)\big)\cr
&<1-{1\over2}a_3\Delta_T,&(2.9)
}
$$
as we may assume that $a_3\Delta_T$ is small 
while $a_3$ is large. Hence we have proved that if $\beta$ exists,
then
$$
\sum_{\scr{q\le Q}\atop\scr{(q,k)=1}}
\mathop{{\sum}^*}_{\chi\bmod q}
\bigg|\sum_{\scr{p\le x}\atop\scr{p\equiv\ell\bmod k}}
\chi(n)\Lambda(n)\bigg|\le 2{x\over\varphi(k)}\Big
(1-{1\over5}a_3\Delta_T\Big),\eqno(2.10)
$$
provided $a_3\Delta_T$ is small and 
$x\ge T^c$ with a computable absolute constant $c>0$.
The case where the $T$-exceptional zero does not exist
is analogous; in fact, simpler. The rest
of the proof may be skipped.
\bigskip
\noindent
{\bf 3.} Both the bounds $(1.2)$ and $(2.5)$ and thus $(1.6)$
are sieve results; 
that is, they are proved using mainly sieve arguments,
without the zero-density theory 
or the Deuring--Heilbronn phenomenon. The proof in [12]
of the assertion (1.2) depends on Iwaniec's work [2] on the
bilinear structure in the error term of the combinatorial linear sieve;
an alternative approach to his
result itself can be found in [12] (see also [1]).
Prior to [2], a bilinear structure in the error term of the Selberg sieve
was observed in [6] and the first uniform
improvement upon $(1.4)$ was 
achieved; see [3]. Later
the development [13] made it possible to 
prove $(1.2)$
via the Selberg sieve as well (see [17] for 
a further development). On the other
hand, the bound $(2.5)$ depends on our large sieve extension [8]
of the Selberg sieve that is devised in [7]
via the duality principle and the quasi-character
derived from optimal $\Lambda^2$-weights.  
This line of consideration yielded a new 
way [10] to discuss zero-free regions of
the zeta-function; in fact, it gave an
assertion [14] that appears beyond the reach of the
convexity argument of Borel, Carath\'eodory and Landau.
A historical account of the developments in the modern
theory of sieve methods can be found in [15][16].
\bigskip
\noindent
{\bf References}
\medskip
\noindent
\item{[1]} G. Greaves. {\it Sieves in Number Theory\/}. 
Springer-Verlag, Berlin 2001.
\item{[2]} H. Iwaniec. A new form of the error 
term in the linear sieve. Acta Arith.,  {\bf37} (1980), 307--320.
\item{[3]} ---. Sieve methods. In:
{\it Intern.\ Congress of Math.\ Proc., Helsinki 1978\/},  
Acad.\ Sci.\ Fennica, Helsinki 1980, pp.\ 357--364.
\item{[4]} J. Maynard. On the Brun--Titchmarsh theorem.
arXiv:1201.1777 (math.NT).
\item{[5]} H.L. Montgomery and R.C. Vaughan. 
The large sieve. Mathematika, {\bf20} (1973), 119--134.
\item{[6]} Y. Motohashi. On some improvements of the 
Brun--Titchmarsh theorem. J. Math.\ Soc.\ Japan, 
{\bf 26} (1974), 306--323.
\item{[7]} ---. On the Deuring--Heilbronn phenomenon. Part I.  
Proc.\ Japan Acad., {\bf53} (1977), 1--2; Part II.  ibid, 25--27.
\item{[8]} ---. A note on the large sieve.  
Proc.\ Japan Acad., {\bf53} (1977), 17--19; Part
II. ibid, 122--124: Part III.  ibid, {\bf55}A
(1979), 92--94: Part IV. ibid,  {\bf56}A
(1980), 288--290.
\item {[9]} ---. A note on Siegel's zeros.  
Proc.\ Japan Acad., {\bf55}A (1979), 190--192.
\item {[10]} ---. An elementary proof of Vinogradov's 
zero-free region for the Riemann zeta-function. 
In: {\it Recent Progress in Analytic Number Theory\/}, Vol.1, 
Academic Press, London 1981, pp.\ 257--267.
\item{[11]} ---. Primes in arithmetic progressions.
Invent.\ math., {\bf 44} (1978), 163--178.
\item{[12]} ---. {\it Sieve Methods and 
Prime Number Theory\/}. LN in Math.\ Phys., {\bf72},
Tata IFR, Bombay 1983. A web edition:
www.math.tifr.res.in/$\sim$publ/ln/
\item{[13]} ---. On the error term in the Selberg sieve.  
In: {\it Number Theory in Progress: A. Schinzel Festschrift\/}, 
Walter de Gruyter,  Berlin, 1999, pp.\ 1053--1064.
\item{[14]} ---. An observation on the 
zero-free region of the Riemann
zeta-function. {\it A. S\'ark\"ozy Festschrift}, 
Periodica Math.\ Hungarica, 
{\bf 42} (2001), 117--122. 
\item{[15]} ---. An overview of sieve
methods and their history.\ Sugaku Expositions.
Amer.\ Math. Soc., {\bf21} (2008),
1--32. arXiv:math/0505521
\item{[16]} ---. {\it Analytic Number Theory\/}.\ I. 
{\it Distribution of Prime Numbers\/}. Asakura Books, 
Tokyo 2009; II. {\it Zeta Analysis\/}. ibid, 2011. (Japanese)
\item{[17]} Y. Motohashi and J. Pintz.
A smoothed GPY sieve.
Bulletin of London Math.\ Soc., {\bf40} (2008), 
pp.\ 298--310. 

\vskip 1cm
\noindent
\font\small=cmr8
{\small\noindent 
Department of Mathematics, Nihon University,
Surugadai, Tokyo 101-8308, Japan
}
\hfill\def\ymzeta
{\font\brm=cmr17 at 30pt\font\sssrm=cmr5 at 4pt
\font\ssssrm=cmr5 at2.5pt
{{\brm O}\raise 9pt\hbox{\hskip -22pt
$\hfil\raise3pt\hbox{\ssssrm KH}\atop\hbox{
{\sssrm Y}$\zeta$\hskip-1pt{\sssrm M}}$}}}
\ymzeta

\bye